\documentclass[11pt,a4paper,reqno]{amsart}
\usepackage{amsmath}
\usepackage{amsfonts}
\usepackage{amssymb}
\usepackage{amscd}
\usepackage{url}
\usepackage[pdftex,bookmarks=true]{hyperref}

\newcommand\Vol{{\operatorname{Vol}}}

\newcommand\R{{\mathbf{R}}}
\newcommand\C{{\mathbf{C}}}

\renewcommand\P{{\mathbf{P}}}
\newcommand\E{{\mathbf{E}}}

\newcommand\Z{{\mathbf{Z}}}

\newcommand\F{{\mathbf{F}}}
\newcommand\I{{\mathbf{I}}}

\newcommand\eps{\varepsilon}

\newcommand\ep{\epsilon}

\theoremstyle{plain}
  \newtheorem{theorem}[subsection]{Theorem}

  \newtheorem{lemma}[subsection]{Lemma}
  \newtheorem{corollary}[subsection]{Corollary}
  
  \newtheorem{example}[subsection]{Example}

\theoremstyle{remark}
  \newtheorem{remark}[subsection]{Remark}

\theoremstyle{definition}
  \newtheorem{definition}[subsection]{Definition}

\include{psfig}

\begin{document}

\title[Inverse Littlewood-Offord Problem]{Optimal Inverse Littlewood-Offord Theorems}

\author{Hoi Nguyen}
\address{Department of Mathematics, Rutgers University, Piscataway, NJ 08854}
\email{hoi@math.rutgers.edu}

\author{Van Vu}
\address{Department of Mathematics, Rutgers University, Piscataway, NJ 08854}
\email{vanvu@math.rutgers.edu}
\thanks{Both authors are supported by research grants DMS-0901216 and AFOSAR-FA-9550-09-1-0167.}
\subjclass[2000]{11B25}
\keywords{inverse Littlewood-Offord problem, concentration probability, generalized arithmetic progression}

\maketitle

\begin{abstract}
Let $\eta_i, i=1,\dots, n$ be iid Bernoulli random variables, taking values $\pm 1$ with probability $\frac{1}{2}$. Given
a multiset $V$ of $n$ integers  $v_1, \dots, v_n$,  we define the \emph{concentration probability} as

$$\rho(V) := \sup_{x} \P( v_1 \eta_1+ \dots
v_n \eta_n=x). $$

A classical result of Littlewood-Offord and Erd\H{o}s 
from the 1940s asserts that, if the $v_i $ are non-zero, then $\rho(V)$ is $O(n^{-1/2})$.  Since then, many researchers have obtained improved bounds by assuming various extra restrictions on $V$.

About 5 years ago, motivated by problems concerning random matrices, Tao and Vu introduced the Inverse Littlewood-Offord problem.
In the inverse problem, one would like to characterize the set $V$, given that $\rho(V)$ is relatively large.

In this paper, we introduce a new method to attack the inverse problem. As an application, we strengthen the previous result of Tao and Vu, obtaining  an optimal characterization for $V$. This immediately implies several classical theorems, such as those of S\'ark\"ozy-Szemer\'edi and Hal\'asz.

The method also applies to the continuous setting and leads to a simple proof for the $\beta$-net theorem of Tao and Vu, which plays a key role in their recent studies of random matrices.

All results extend to the general case when $V$ is a subset of an abelian torsion-free group, and $\eta_i$ are independent variables satisfying some weak conditions.

\end{abstract}

\section{Introduction}\label{section:introduction}

\subsection{The Forward Littlewood-Offord problem}

Let $\eta_i, i=1,\dots, n$ be iid Bernoulli random variables, taking values $\pm 1$ with probability $\frac{1}{2} $. Given
a multiset $V$ of $n$ integers  $v_1, \dots, v_n$, we define the random
walk $S$ with steps in $V$ to be the random variable $S:=\sum_{i=1}^n{v_i\eta_i}$. The \emph{concentration probability} is defined to be

$$\rho(V) := \sup_{x} \P(S=x). $$

Motivated by their study of random polynomials in the 1940s, Littlewood and
Offord \cite{LO} raised the question of bounding $\rho(V)$. (We call this the  {\it forward} Littlewood-Offord problem, in contrast with the {\it inverse} Littlewood-Offord problem discussed in the next section.) They
showed that $\rho(V)=O(n^{-1/2}\log n)$.
Shortly after the Littlewood-Offord paper, Erd\H{o}s \cite{ELO}
gave a beautiful combinatorial proof of the
refinement

\begin{equation} \label{eqn:Erdos} 
\rho(V)\le \frac{ \binom{n}{n/2}}{2^n } =O(n^{-1/2}). 
\end{equation}

Erd\H os' result is sharp, as demonstrated by  $V=\{1,\dots,1\}$.

\vskip2mm

{\it Notation.} Here and later, asymptotic notations, such as $O, \Omega, \Theta$, and so forth, are used under the assumption that $n \rightarrow \infty$. A notation such as $O_{C}(.)$ emphasizes that the hidden constant in $O$ depends on $C$. If $a= \Omega (b)$, we write $b \ll a$ or $a \gg b$.  All logarithms have a natural base, if not specified otherwise.

\vskip2mm

The results of Littlewood-Offord and Erd\H{o}s are
classics in combinatorics and have generated an impressive wave of research, 
particularly from the early 1960s to the late 1980s.

One direction of research was to generalize Erd\H os' result to
other groups. For example, in 1966 and 1970, Kleitman extended Erd\H{o}s' result to complex numbers and normed vectors, respectively. Several results in this direction can be found in \cite{Griggs,Kat}.

Another direction was motivated by the observation that \eqref{eqn:Erdos}  can be improved significantly by making additional assumptions about $V$. The first such result was discovered by Erd\H{o}s and Moser \cite{EM}, who showed that if
$v_i$ are distinct, then $\rho(V) = O(n^{-3/2} \log n)$.
They conjectured that the logarithmic term is not necessary, and this
was confirmed by S\'ark\"ozy and Szemer\'edi
\cite{SS}.

\begin{theorem} \label{theorem:SS} Let $V$ be a set of $n$ different integers, then

$$\rho(V)=O(n^{-3/2}).$$ \end{theorem}

In \cite{H} (see also in \cite{TVbook}), Hal\'asz
proved very general theorems that imply Theorem \ref{theorem:SS} and many others.
One of his results can be formulated as follows.

\begin{theorem} \label{theorem:Hal}
Let $l$ be a fixed integer and $R_l$ be the number of solutions of the equation
$v_{i_1}+\dots+v_{i_l}=v_{j_1}+\dots+v_{j_l}$. Then

$$\rho(V)= O(n^{-2l-\frac{1}{2}}R_l).$$
\end{theorem}

It is easy to see, by setting $l=1$, that Theorem \ref{theorem:Hal} implies Theorem \ref{theorem:SS}.

Another famous result in this area is that of Stanley \cite{Stan}, which, solving a conjecture of Erd\H os and Moser, shows when $\rho(V)$ attains its maximum under the assumption that the $v_i$ are different.

\begin{theorem} \label{theorem:Stan}
Let $n$ be odd and $V_0 :=\{- \lfloor n/2 \rfloor , \dots, \lfloor n/2 \rfloor \}$. Then $$\rho(V) \le \rho(V_0) . $$
\end{theorem}

A similar result holds for the case of $n$ being even \cite{Stan}.
Stanley's proof of Theorem \ref{theorem:Stan} used sophisticated machinery from algebraic geometry, particularly the hard Lefschetz theorem.
A few years later, a more elementary proof was given by Proctor \cite{P}. This proof also has an algebraic nature, involving the representation of the Lie algebra $sl(2, \C)$. As far as we know, there is no purely combinatorial proof.

It is natural to ask for the actual value of
$ \rho(V_0)$. From Theorem \ref{theorem:SS}, one would guess (under the assumption that the elements of $V$ are different) that

$$  \rho(V_0) = (C_0+o(1)) n^{-3/2}  $$ for some constant $C_0 >0$. However, the algebraic proofs do not give the value of $C_0$. In fact, it is not obvious that $\lim_{n \rightarrow \infty} n^{3/2} \rho(V_0)$ exists.

Assuming that $C_0$ exists for a moment, one would next wonder if $V_0$ is a stable maximizer. In other words, if some other set $V_0'$ has $\rho(V_0')$  close to $C_0 n^{-3/2}$, then should $V_0'$ (possibly after a normalization) be ''close'' to $V_0 $ ?
(Note that $\rho$ is invariant under dilation, so a normalization would be necessary.)

\subsection{The inverse Littlewood-Offord problem} Motivated by inverse theorems from additive combinatorics (see \cite[Chapter 5]{TVbook}) and a variant for random sums in \cite[Theorem 5.2]{TVsing}, Tao and the second author \cite{TVinverse} brought a different view to the problem. Instead of trying to improve the bound further by imposing new assumptions (as done in the forward problems), they tried to provide the complete picture by finding the underlying reason as to why the concentration probability is large (say, polynomial in $n$).

Note that the (multi)-set  $V$ has $2^n$ subsums, and $\rho({V})\ge n^{-C}$ means that at least $\frac{2^n}{n^C}$ of these take the same value. This observation suggests that the set should have a very strong additive structure. To determine this structure,
we first discuss a few examples of $V$, where $\rho(V)$ is large. For a
set $A$, we denote the set $ \{a_1+ \dots + a_l| a_i \in A\}$ by $lA$.

\begin{example}  Let $I=[-N,N]$ and $v_1, \dots, v_n$ be elements of
$I$. Because $S \in nI$, by the pigeon-hole principle, $\rho(V)  \ge
\frac{1}{|nI|} = \Omega (\frac{1}{nN})$.  In fact, a short
consideration yields a better bound. Note that, with a probability of
least $.99$, we have $S \in 10 \sqrt n I$. Thus, again by the pigeon-hole
principle, we have $\rho(V)  = \Omega (\frac{1}{\sqrt n N})$. If we
set $N=n^{C-1/2}$ for some constant $C\ge 1/2$, then
\begin{equation}\label{bound1} 
\rho(V)  = \Omega (\frac{1}{n^{C}}).
\end{equation}
\end{example}

The next, and more general, construction comes from additive
combinatorics. A very important concept in this area is that of  \emph{generalized 
arithmetic progressions} (GAPs). A set $Q$ is a \emph{GAP of 
rank $r$} if it can be expressed as in the form
$$Q= \{a_0+ x_1a_1 + \dots +x_r a_r| M_i \le x_i \le M_i' \hbox{ for all } 1 \leq i \leq r\}$$
for some $\{a_0,\ldots,a_r\}, \{M_1,\ldots,M_r\},$ and $\{M'_1,\ldots,M'_r\}$.

It is convenient to think of $Q$ as the image of an integer box $B:= \{(x_1, \dots, x_r) \in \Z^r| M_i \le m_i\le M_i' \} $ under the linear map
$$\Phi: (x_1,\dots, x_r) \mapsto a_0+ x_1a_1 + \dots + x_r a_r. $$
The numbers $a_i$ are the \emph{generators } of $P$, the numbers $M_i$ and $M_i'$ are the \emph{dimensions} of $P$, and $\Vol(Q) := |B|$ is the \emph{volume} of $B$. We say that $Q$ is \emph{proper} if this map is one-to-one or, equivalently, if $|Q| = \Vol(Q)$.  For non-proper GAPs, we, of course, have $|Q| < \Vol(Q)$.
If $-M_i=M_i'$ for all $i\ge 1$ and $a_0=0$, we say that $Q$ is {\it symmetric}.

\begin{example} \label{mainexample} Let $Q$ be a proper symmetric GAP of rank $r$ and volume $N$.
Let $v_1, \dots, v_n$ be (not necessarily distinct) elements of
$P$. The random variable $S =\sum_{i=1}^n v_i\eta_i$ takes values in
the GAP $nP$. Because $|nP| \le \Vol (nB) = n^r N$, the pigeon-hole
principle implies that $\rho(V)   \ge \Omega (\frac{1}{n^r N})$. In
fact, using the same idea as in the previous example, one can
improve the bound to $\Omega (\frac{1}{n^{r/2} N})$. If we set
$N=n^{C-r/2}$ for some constant $C\ge r/2$, then

\begin{equation}\label{bound2} 
\rho(V)   = \Omega (\frac{1}{n^{C}}).
\end{equation}
\end{example}

The examples above show that, if the elements of $V$
belong to a proper GAP with a small rank and small cardinality, then
$\rho(V)$ is large. A few years ago, Tao and the second author  \cite{TVinverse}
showed that this is essentially the only reason:

\begin{theorem}[Weak inverse theorem]\label{theorem:weak} \cite{TVinverse} Let  $C, \epsilon > 0$ be arbitrary constants.
There are constants $r$ and $C'$ depending on $C$ and $\epsilon$
such that the following holds.
 Assume that $V  = \{v_1, \ldots, v_n\}$ is a multiset of integers satisfying
$\rho(V) \geq n^{-C}$. Then, there is a proper symmetric GAP $Q$ with a rank of at most $r$ and a volume of at most $n^{C'}$ that contains all but at most
$n^{1-\epsilon}$ elements of $V$ (counting multiplicity).

\end{theorem}

\begin{remark}
The presence of a small set of exceptional elements is 
not completely avoidable. For instance, one can add $o(\log n)$ completely
arbitrary elements to $V$ and, at worst, only decrease $\rho(V) $ by a factor 
of $n^{-o(1)}$.  Nonetheless, we expect the number of such elements to be less than what is given by the results here.
\end{remark}

The reason we call Theorem \ref{theorem:weak} {\it weak} is
that $C'$ is not optimal. In particular, it is far from reflecting
the relations in \eqref{bound1} and \eqref{bound2}.
In a later paper \cite{TVstrong}, Tao and the second author refined  the approach to obtain the following stronger result.

\begin{theorem}[Strong inverse theorem] \label{theorem:strong}
\cite{TVstrong}  Let $C$ and $1> \eps$ be positive constants.
Assume that
$$\rho (V)  \ge n^{-C}. $$
 Then, there exists a proper symmetric  GAP $Q$ of rank $r= O_{C, \eps} (1)$ that contains all but $O_r(n^{1 -\eps} )$
 elements of $V$ (counting multiplicity), where
 $$|Q| = O_{C, \eps} (n^{C - \frac{r}{2} + \eps}). $$
\end{theorem}

The bound on $|Q|$ matches Example \ref{mainexample}, up to the $n^{\epsilon}$
term. However, this error term seems to be the limit of the approach.
The proofs of Theorems \ref{theorem:weak} and
\ref{theorem:strong} rely on a replacement argument and  various lemmas about random walks 
and GAPs.

Let us now consider an application of Theorem 
\ref{theorem:strong}. Note that  Theorem \ref{theorem:strong}
 enables us to make very precise counting arguments. Assume that we would like to count the number of (multi)sets $V$ of integers with $\max |v_{i}| \le N=n^{O(1)}$ such that $\rho(V) \ge \rho:= n^{-C}$.

Fix $d \ge 1$, and fix \footnote{A more detailed version of Theorems
 \ref{theorem:weak} and \ref{theorem:strong} tells us that there are
 not too many ways to choose the generators of $Q$. In particular, if $N =n^{O(1)}$, the number of ways to fix these is 
negligible compared to the main term.} a  GAP $Q$  with rank $r$ and volume $|Q| = n^{C -\frac{r}{2}}$. The dominating term in the calculation will be the number of multi-subsets of size $n$  of $Q$, which is

\begin{equation}\label{discretcounting}  
|Q|^{n }= n^{(C-\frac{r}{2} +\epsilon)n} \le n^{Cn} n^{-\frac{n}{2}+\epsilon n}= \rho^{-n} n^{-n(\frac{1}{2}-\epsilon) }.
\end{equation}

Motivated by questions from random matrix theory, Tao and the second author
obtained the following continuous analogue of this result. 

\begin{definition}[Small ball probability] Let $z$ be a real random variable, and let $V=\{v_1,\dots,v_n\}$ be a multiset in $\R^d$. For any $r>0$, we define the {\it small ball probability} as

$$\rho_{r,z}(V):=\sup_{x\in \R^d} \P(v_1z_1+\dots v_nz_n \in B(x,r)),$$

where $z_1,\dots, z_n$ are iid copies of $z$, and $B(x,r)$ denotes the closed disk of radius $r$ centered at $x$ in $\R^d$. 

\end{definition}

Let $n$ be a positive integer and $\beta,\rho$ be positive numbers that may depend on $n$. Let $\mathcal{S}_{n,\beta, \rho}$ be the collection of all multisets $V=\{v_1,\dots,v_n\} , v_i \in \R^2 $ such that  $\sum_{i=1}^n \| v_i\| ^2=1$ and $\rho_{\beta,\eta}(V)\ge \rho$, where $\eta$ has a Bernoulli distribution. 

\vskip2mm

\begin{theorem}[The $\beta$-net Theorem] \cite{TVcir}\label{theorem:continuous:TV}
Let  $0 <\ep \le 1/3$ and $C>0$ be constants. Then, for all sufficiently large $n$  and $\beta \ge \exp(-n^{\ep})$
and $\rho \ge n^{-C}$, there is a set $\mathcal{S}\subset (\R^2)^n$ of size at most

$$\rho^{-n}n^{-n(\frac{1}{2}-\epsilon) }
  + \exp(o(n))$$

\noindent such that for any $V=\{v_1,\dots,v_n\} \in \mathcal{S}_{n,\beta,\rho}$, there is some $V'=(v_1',\dots, v_n')\in \mathcal{S}$ such that $\|v_i-v_i'\|_{2}  \le \beta$ for all $i$.
\end{theorem}

The theorem looks a bit cleaner if we use $\C$ instead of  $\R^2$ (as in \cite{TVcir}). However, we prefer the current form, because it is more suitable for generalization. 
The set   $\mathcal{S}$ is usually referred to as a $\beta$-net of  $ \mathcal{S}_{n,\beta,\rho}$.

Theorem \ref{theorem:continuous:TV} is at the heart of establishing the Circular Law conjecture in random matrix theory (see \cite{TVcir, TVbull}). It also plays 
an important role in the study of the condition number of randomly perturbed matrices (see \cite{TVcomp}). 
Its proof in \cite{TVcir} is quite technical and occupies the bulk of that paper.

However, given the above discussion, one might expect to obtain Theorem \ref{theorem:continuous:TV}   as a simple corollary of
a continuous analogue of Theorem \ref{theorem:strong}.  However, the arguments in \cite{TVcir} have not yet provided such an inverse theorem
(although they  did provide a sufficient amount of information about the set $S$ to make an estimate possible).  The paper \cite{RV} by Rudelson and Vershynin also contains a characterization of the set $S$, but their characterization has a somewhat different spirit than those discussed in this paper.

\section{A new approach and new results}

In this paper, we introduce a new approach to the inverse theorem. The core of this new approach is a (long-range) variant of Freiman's famous inverse theorem.

This new approach seems powerful.  First, it enables us to remove the error term $n^{\epsilon}$ in Theorem \ref{theorem:strong}, resulting in an optimal inverse theorem.

\begin{theorem}[Optimal inverse Littlewood-Offord theorem, discrete case]\label{theorem:optimal}
Let $\eps<1$ and $C$ be positive constants. Assume that
$$\rho (V)  \ge  n^{-C}. $$ Then, there exists a proper symmetric GAP $Q$ of rank $r= O_{C, \eps} (1)$ that contains all but at most $\eps n$
 elements of $V$ (counting multiplicity), where
 $$|Q| = O_{C, \eps} ( \rho (V)^{-1} n^{- \frac{r}{2}}). $$
\end{theorem}

This immediately implies several forward theorems, such as Theorems \ref{theorem:SS} and \ref{theorem:Hal}. For example, we can prove Theorem \ref{theorem:SS} as follows.

\begin{proof}(of Theorem \ref{theorem:SS}) Assume, for contradiction, that there is a set $V$ of 
$n$ distinct numbers such that $\rho(V) \ge c_1 n^{-3/2}$ for some large constant $c_1$ to be chosen. Set $\eps=.1, C=3/2$.
By Theorem \ref{theorem:optimal}, there is a GAP $Q$ of rank $r$ and size $O_{C, \epsilon} (\frac{1}{c_1} n^{C-\frac{r}{2}} )$ that contains
at least $.9n$ elements from $V$. This implies $|Q| \ge .9 n$. By setting $c_1$ to be sufficiently large and using the fact that   $C=3/2$ and $r \ge 1$, we can guarantee that $|Q|  \le .8n $, a contradiction.
\end{proof}

Theorem \ref{theorem:Hal} can be proved in a similar manner with the details left as an exercise.

Similar to \cite{TVstrong, TVinverse}, our method and results can be extended (rather automatically) 
to much more general settings.

{\it General $V$.} Instead of taking $V$ to be a subset of $\Z$, we can take it to be a subset of any abelian torsion-free group $G$
(thanks to Freiman isomorphism, see Section \ref{section:GAPs}). We can also replace $\Z$ by the finite field $\F_p$, where $p$ is any sufficiently large prime. (In fact, the first step in our proof is to embed
$V$ into $\F_p$.)

{\it General $\eta$.} We can replace the Bernoulli random variables by independent random variables $\eta_i$ satisfying the following condition. There is a constant $c >0$
and an infinite sequence of primes $p$ such that for any $p$ in the sequence, any (multi)-subset
$V$ of size $n$ of $\F_p$  and any $t \in \F_p$

\begin{equation} \label{eqn:condition0} \prod_{i=1}^n |\E e_p(\eta_i  v_i  t )| \le \exp(- c \sum_{i=1}^n \| \frac{v_i t}{p} \| ^2 ) \end{equation}  where $\|x\| $ denotes the distance from $x$ to the closest integer
(we view the elements of $\F_p$ as integers between $0$ and $p-1$) and $e_p(x) := \exp(2\pi \sqrt{-1} x/p)$.

\begin{example} (Lazy random walks) Given a parameter  $0 <  \mu \le 1$,
let  $\eta_i^\mu$ be  iid copies of a random
variable $\eta^\mu$, where $\eta^\mu=1$ or $-1$ with probability $\mu/2$,
and $\eta^\mu=0$ with probability $1-\mu$.  The sum

$$S^\mu(V):=\sum_{i=1}^n {\eta_i}^\mu v_i,$$

\noindent  can be viewed as a {\it lazy random walk} with steps in $V$.  A simple calculation shows

$$\E e_p(\eta x)=  (1-\mu) + \mu \cos \frac{2 \pi x}{p}  . $$

It is easy to show that there is a constant $c >0$ depending on $\mu$ such that

$$| (1-\mu) + \mu \cos \frac{2 \pi x}{p}  |  \le \exp( -c \|\frac{x}{p} \|^2) . $$
\end{example}

\begin{example} ($\mu$-bounded variables) It suffices to assume that there is some constant  $0 < \mu \le 1$ such that for all $i$

\begin{equation}\label{eqn:condition} 
|\E e_p(\eta_i  x) | \le (1-\mu) + \mu \cos  \frac{2 \pi x}{p}. 
\end{equation}

\end{example}

\begin{theorem} \label{theorem:general}  The conclusion of Theorem \ref{theorem:optimal} holds 
for the case when $V$ is a multi-subset of an arbitrary  torsion-free abelian group $G$ and $\eta_i, 1\le i \le n$ are  independent random variables satisfying \eqref{eqn:condition0}.
\end{theorem}

In some applications, we might need a version of Theorem \ref{theorem:optimal} with a smaller number of exceptional elements. By slightly modifying the proof presented in Section \ref{section:optimal}, we can prove the following result.

\begin{theorem}\label{theorem:optimal:general} Let $\eps<1$ and $C$ be positive constants. Assume that

$$\rho (V)  \ge  n^{-C}. $$

Then, for any $n^\ep \le n' \le n$, there exists a proper symmetric GAP $Q$ of rank $r=O_{\ep,C}(1)$ that contains all but $n'$ elements of $V$ (counting multiplicity), where 

$$|Q|=O_{C,\ep}(\rho^{-1}/{n'}^{r/2}).$$ 

\end{theorem}

\vskip .2in

\begin{remark}
In an upcoming paper \cite{Ng}, we are able to address the unresolved issues concerning Theorem \ref{theorem:Stan} by following the method used to prove Theorem \ref{theorem:optimal}. We prove that $\rho (V_0) = (\sqrt {\frac{24}{\pi}} +o(1)) n^{-3/2}$. More important, we obtain a stable version of Theorem \ref{theorem:Stan}, which shows that, if $\rho(V) $ is close to $(\sqrt {24 /\pi} +o(1)) n^{-3/2} $, then $V$ is ''close'' to $V_0$. As a byproduct, we obtain the first  {\it non-algebraic} proof for the asymptotic version of the Stanley theorem.
\end{remark}

We now turn to the continuous setting.  In this part, we consider a real random variable  $z$ such that there exists a constant $C_z$ such that

\begin{equation}\label{eqn:2ndmoment}
\P(1\le |z_1-z_2| \le C_z)\ge 1/2,
\end{equation}

\noindent where $z_1,z_2$ are iid copies of $z$. We note that Bernoulli random variables are clearly of this type. (Also, the interested reader may find \eqref{eqn:2ndmoment} more general than the condition of the $\kappa$-controlled second moment defined in \cite{TVcir} and the condition of bounded third moment in \cite{RV}.) In the statement above, $C_z$ is not uniquely defined. In what follows, we will 
take the smallest value of  $C_z$.

We say that a vector $v \in   \R^d$ is {\it $\delta$-close} to a set $Q\subset \R^d$ if there exists a vector $q\in Q$ such that $\|v-q\|_2 \le \delta$.
A set $X$ is  $\delta$-close to a set $Q$ if every element of $X$ is $\delta$-close to $Q$.
The analogue of Example \ref{mainexample} is the following.

\begin{example} \label{mainexample2}
Let $Q$ be a proper symmetric GAP of rank $r$ and volume $N$ in $\R^d$.
Let $v_1, \dots, v_n$ be (not necessarily distinct) vectors that are
$O(\beta  n^{-1/2})$-close to $Q$.  If we set
$|Q|=n^{C-\frac{r}{2}}$ for some constant $C\ge r/2$, then
\begin{equation}\label{bound3} 
\rho_{\beta, \eta} (V)   = \Omega (\frac{1}{n^{C}}).
\end{equation}
\end{example}

Thus, one would expect that, if $\rho_{\beta, z} (V)   $ is large, then (most of) 
$V$ is $O(\beta n^{-1/2})$-close to a GAP with a small volume.
Confirming this intuition, we obtain the following continuous analogue of
Theorem \ref{theorem:optimal}.


\begin{theorem}[Optimal inverse Littlewood-Offord theorem, continuous  case]\label{theorem:optimal1} Let $\delta, C>0$ be arbitrary
constants and $\beta>0 $ be a parameter that may depend on $n$.
 Suppose that $V=\{v_1,\dots,v_n\}$ is a (multi-) subset of $\R^d$ such that $\sum_{i=1}^n\|v_i\|_2^2=1$ and that $V$ has large small ball probability
$$\rho:= \rho_{\beta,z}(V)\ge n^{-C}, $$ where $z$ is a real  random variable satisfying \eqref{eqn:2ndmoment}.
 Then, there exists a proper symmetric GAP $Q$ of rank $d\le r= O(1)$ so that 
all but at most $\delta n$ elements of $V$ (counting multiplicity) are $O(\beta \frac{\log n}{n^{1/2}}) $-close to $Q$, where

 $$|Q| = O ( \rho ^{-1}  \delta^{(-r+d)/2} n^{(-r+d)/2} ). $$

\end{theorem}

The theorem is optimal in the sense that the exponent $(-r+d)/2$ of $n$ cannot generally be improved (see Appendix \ref{section:corcont}  for more details).

Theorem \ref{theorem:optimal1} is a special case of the following more general theorem.

\begin{theorem}[Continuous Inverse Littlewood-Offord theorem, general setting]\label{theorem:continuous:NgV} Let $0 <\ep < 1; 0 < C$ be constants. Let 
$ \beta >0$ be a parameter that may depend on $n$. Suppose that $V=\{v_1,\dots,v_n\}$ is a (multi-) subset of $\R^d$ such that $\sum_{i=1}^n\|v_i\|_2^2=1$ and that $V$ has large small ball probability
$$\rho:= \rho_{\beta,z}(V)\ge n^{-C}, $$ where $z$ is a real  random variable satisfying \eqref{eqn:2ndmoment}. Then, the following holds. For any number $n^\ep \le n'  \le n$, there exists a proper symmetric GAP $Q=\{\sum_{i=1}^r x_ig_i : |x_i|\le L_i \}$ such that

\begin{itemize}

\item (Full dimension) There exists $\sqrt{\frac{n'}{\log n}}\ll k \ll \sqrt{n'}$ such that the dilate $P:= \beta^{-1} k \cdot Q$ contains the discrete hypercube $\{0,1\}^d $.

\vskip .1in

\item (Approximation) At least $n-n'$ elements of $V$ are $O(\frac{\beta}{k})$-close to $Q$.

\vskip .1in

\item (Small rank and cardinality) $Q$ has constant rank $d \le r=O (1)$, and cardinality

$$|Q| =O(\rho^{-1} n'^{(-r+d)/2}).$$

\vskip .1in

\item (Small generators) There is a non-zero integer $p=O(\sqrt{n'})$ such that  all steps $g_i$ of $Q$ have the form  $g_i=(g_{i1},\dots,g_{id})$, where $g_{ij}=\beta \frac{p_{ij}} {p} $ with  $p_{ij} \in \Z$ and $p_{ij}=O(\beta^{-1} \sqrt {n'}  ).$

\end{itemize}
\end{theorem}

Theorem \ref{theorem:continuous:NgV}  implies the following corollary (see Appendix \ref{section:corcont} for a simple proof), from which one can derive  Theorem \ref{theorem:continuous:TV} in a straightforward manner (similar to the discrete case discussed earlier).

\begin{corollary}\label{corollary:continuous:NgV} Let $0 <\ep < 1; 0 < C$ be constants. Let 
$ \beta >0$ be a parameter that may depend on $n$. Suppose that $V=\{v_1,\dots,v_n\}$ is a (multi-) subset of $\R^d$ such that $\sum_{i=1}^n\|v_i\|_2^2=1$ and that $V$ has large small ball probability
$$\rho:= \rho_{\beta,z}(V)\ge n^{-C}, $$ where $z$ is a real  random variable satisfying \eqref{eqn:2ndmoment}. Then the following holds. For any number $n'$ between $n^\ep$ and $n$, there exists a proper symmetric GAP $Q=\{\sum_{i=1}^r x_ig_i : |x_i|\le L_i \}$ such that

\begin{itemize}

\item At least $n-n'$ elements of $V$ are $\beta$-close to $Q$.

\vskip .1in

\item $Q$ has small rank, $r=O (1)$, and small cardinality

$$|Q| \le \max \left(O(\frac{\rho^{-1}}{\sqrt{n'}}),1\right).$$

\vskip .1in

\item There is a non-zero integer $p=O(\sqrt{n'})$ such that all
 steps $g_i$ of $Q$ have the form  $g_i=(g_{i1},\dots,g_{id})$, where $g_{ij}=\beta \frac{p_{ij}} {p} $ with $p_{ij} \in \Z$ and $p_{ij}=O(\beta^{-1}\sqrt{n'}).$

\end{itemize}
\end{corollary}

Note that the approximations obtained from Corollary \ref{corollary:continuous:NgV} are rougher than those from Theorem \ref{theorem:continuous:NgV}). However, the bound on $|Q|$ is improved in some critical cases (particularly when $r=d$).

In the above theorems, the hidden constants could depend on previously set constants $\epsilon, C, C_z, d$.  We could have written
$O_{\epsilon, C ,C_z, d}$ and $\ll_{\epsilon, C ,C_z, d}$ everywhere, but these notations are somewhat cumbersome, and this dependence is not our focus.

\begin{proof}(of Theorem \ref{theorem:continuous:TV}) Set $n':=n^{1-\frac{3\ep}{2}}$ (which is $\gg n^{\ep}$ as $\ep\le 1/3$). Let $\mathcal{S'}$ be the collection of
 all subsets of size at least $n-n'$ of GAPs whose parameters satisfy the conclusion of Corollary \ref{corollary:continuous:NgV}.

Because each GAP is determined by its generators and dimensions, the number of such GAPs is bounded by $((\beta^{-1}\sqrt{n'})\sqrt{n'})^{O(1)}
(\frac{\rho^{-1}}{\sqrt{n'}})^{O(1)}=\exp(o(n))$.  (The term $(\frac{\rho^{-1}}{\sqrt{n'}})^{O(1)}$
bounds the number of choices of the dimensions $M_i$.)  Thus, $|\mathcal{S'}|= \left(O((\frac{\rho^{-1}}
{\sqrt{n'}})^{n})+1\right) \exp(o(n))$.

We approximate each of the exceptional elements by a lattice point in $\beta\cdot (\Z/d)^d$.
Thus, if we let $\mathcal{S''}$ to be the set of these approximated tuples, then $|\mathcal{S''}|\le \sum_{i\le n'} (O(\beta^{-1}))^i = \exp(o(n))$ (here, we used the assumption $\beta \ge \exp(-n^{\ep})$).

Set $\mathcal{S}:=\mathcal{S'}\times \mathcal{S''}$. It is easy to see that
 $|\mathcal{S}|\le O(n^{-1/2+\ep}\rho^{-1})^n+\exp(o(n))$. Furthermore, if $\rho(V)\ge n^{-O(1)}$, then $V$ is $\beta$-close to an element of $\mathcal{S}$, concluding the proof.
\end{proof}

\section{The long range inverse theorem}

Let us first recall a famous theorem by Freiman \cite[Chapter 5]{TVbook}.

\begin{theorem}[Freiman's inverse theorem]\label{theorem:Freiman}
Let $\gamma$ be a positive constant and  $X$ a subset of a torsion-free group such that $|2X| \le \gamma |X|$. Then, there
is a proper symmetric GAP $Q$ of rank at most $r=O_\gamma(1) $ and cardinality $O_\gamma(|X|)$ such that $ X \subset Q$.
\end{theorem}

In our analysis, we will need to deal with an assumption of the form $|kX| \le k^\gamma|X|$, where $\gamma$ is a constant but
$k$ is not. (Typically, $k$ will be a positive power of $|X|$.)  We successfully give a structure for $X$ under this condition
in the following theorem, which we will call the long range inverse theorem.

\begin{theorem}[Long range inverse theorem]\label{theorem:longrange}
Let $\gamma>0$ be constant. Assume that $X$ is a subset of a torsion-free group such that $0 \in X$ and $ |kX| \le k^\gamma|X|$ for some integer $k \ge 2$ that may depend on $|X|$. Then, there is proper symmetric GAP $Q$ of rank $r=O(\gamma)$ and cardinality $O_{\gamma}( k^{-r} |kX| )$ such that $X \subset Q$.
\end{theorem}

Note that for any given $\ep>0$ and for any sufficiently large $k$, it is implied from Theorem \ref{theorem:longrange} that the rank of $Q$ is at most $\gamma+\ep$. The implicit constant involved in the size of $Q$ can be taken to be $2^{2^{2^{O(\gamma)}}}$, which is quite poor. Although we have not elaborated on this bound substantially, our method does not seem to say anything when the polynomial growth with a size of $kX$ is replaced by something faster.

Theorem \ref{theorem:longrange} will serve as our main technical tool. This theorem can be proved by applying an earlier result \cite{TVjohn}. We give a short deduction in Appendix \ref{section:longrange}.

\section{Freiman isomorphism}\label{section:GAPs}

We now introduce the concept of Freiman isomorphism that allows us to transfer an additive problem to another group in a way that is more flexible than the usual notion of group isomorphism.

\begin{definition}[Freiman isomorphism of order $k$]
Two sets $V,V'$ of additive groups $G,G'$ (not necessarily
torsion-free) are a Freiman isomorphism of order $k$ (in generalized form) if there is an 
injective map $f$ from $V$ to $V'$ such that $f(v_1)+\dots +f(v_k) = 
f(v_1')+\dots +f(v_k')$ in $G'$ if and only if $v_1+\dots+v_k =
v_1'+\dots +v_k'$ in $G$.
\end{definition}

\noindent The following theorem allows us to pass from an arbitrary
torsion-free group to $\Z$ or cyclic groups of a prime order (see \cite[Lemma 5.25]{TVbook}).

\begin{theorem} Let $V$ be a finite subset of a torsion-free additive group $G$. Then, for
any integer $k$, there is a Freiman isomorphism $\phi$ : $V\rightarrow \phi(V)$ of order $k$ to some
finite subset $\phi(V)$ of the integers $\Z$. The same is true if we replace $\Z$ by $\F_p$, if $p$ is sufficiently large, depending on $V$.
\end{theorem}

An identical proof to that in \cite{TVbook} implies the following stronger result.

\begin{theorem}\label{theorem:Freimaniso}
Let $V$ be a finite subset of a torsion-free additive group $G$. Then, for
any integer $k$, there is a map $\phi$ : $V\rightarrow \phi(V)$ to some
finite subset $\phi(v)$ of the integers $\Z$ such that

\begin{equation}\label{eqn:Freimaniso}
v_1+\dots +v_i = v_1'+\dots + v_j' \Leftrightarrow \phi(v_1)+\dots +\phi(v_i) = \phi(v_1')+\dots \phi(v_j') 
\end{equation}

\noindent for all $i,j\le k$. The same is true if we replace $\Z$ by $\F_p$, if $p$ is sufficiently large, depending on $V$.
\end{theorem}

By Theorem \ref{theorem:Freimaniso}, a large prime $p$ and set $V_p \subset \F_p$ exist such that \eqref{eqn:Freimaniso} holds for all $i,j \le |V|$. Hence, we infer that

$$\rho(V) = \rho(V_p).$$

Thus, instead of working with a subset $V$ of a torsion-free group, it is sufficient to work with a subset of $\F_p$, where $p$ is sufficiently large.

To end this section, we record a useful fact about GAPs, as follows. Assume that $A$ is a dense subset of a GAP $Q$. Then, the iterated sumsets  $kA$ contain a structure similar to $Q$ (see \cite[Lemma 4.4]{SV1}, \cite[Lemma B3]{T-sol}).

\begin{lemma}[S\'ark\"ozy-type theorem in progressions]\label{lemma:Sarkozy} Let
$Q= \{a_1 x_1 + \dots + a_r x_r: |x_i| \le M_i, 1\le i \le r \}$ be a proper GAP in a torsion-free group of rank $r$. Let $A\subset Q$ be a symmetric subset such that
$|A|\ge \delta |Q|$ for some $0 < \delta < 1$. Then, there exists positive integers $1\le m,l \ll_{\delta,r} 1$
such that $Q_l \subset 2mA$, where $Q_l$ is the GAP

$$Q_l = \{la_1 x_1 + \dots + la_r x_r: |x_i| \le M_i/l^2, 1\le i \le r \}.$$
\end{lemma}



\section{Proof of Theorem \ref{theorem:optimal}}\label{section:optimal}

{\it Embedding.} The first step is to embed the problem into the finite field $\F_p$ for some prime $p$. In the case when the $v_i$ are integers, we simply take $p$ to be a large prime 
(for instance, $p \ge 2^n (\sum_{i=1}^n |v_i| +1)$ suffices). If $V$ is a subset of a general torsion-free group $G$, one can use Theorem \ref{theorem:Freimaniso}.

From now on, we can assume that $v_i$ are elements of $\F_p$ for some large prime $p$. 
We view elements of $\F_p$ as integers between $0$ and $p-1$. We use the shorthand $\rho$ to denote $\rho (V)$.

{\it Fourier Analysis.} The main advantage of working in $\F_p$ is that one can use discrete Fourier analysis. Assume that $$\rho= \rho(V)=\P( S=a), $$ for some $a  \in \F_p$.
Using the standard notation $e_p(x)$ for $\exp(2\pi \sqrt{-1} x/p )$, we have

\begin{equation}\label{eqn:fourier1} \rho= \P(S=a)= \E \frac{1}{p} \sum_{\xi\in \F_p} e_p (\xi (S-a)) = \E \frac{1}{p} \sum_{\xi\in \F_p} e_p (\xi S) e_p(-\xi a) .\end{equation}

By independence,

\begin{equation}\label{eqn:fourier2}  
\E e_p(\xi S) = \prod_{i=1}^n e_p(\xi \eta_i v_i)= \prod_{i=1}^n \cos \frac{2\pi \xi v_i}{p}.  
\end{equation}

It follows that

\begin{equation}\label{eqn:fourier3} 
\rho  \le \frac{1}{p} \sum_{\xi \in \F_p} \prod_i |\cos \frac{2 \pi v_i \xi}{p}  |  = \frac{1}{p} \sum_{\xi \in \F_p} \prod_i |\frac{\cos  \pi v_i \xi}{p}|, 
\end{equation}

where we made the variable change $\xi \rightarrow \xi/2$ (in $\F_p$) to obtain the last identity. 

 By convexity, we have that   $|\sin  \pi z | \ge 2 \|z\|$ for any $z\in \R$, where $\|z\|:=\|z\|_{\R/\Z}$ is the distance of $z$ to the nearest integer. Thus,

\begin{equation}\label{eqn:fourier3-1}
|\cos \frac{\pi x}{p}|  \le  1- \frac{1}{2} \sin^2 \frac{\pi x}{p}  \le 1 -2 \|\frac{x}{p} \|^2  \le \exp( - 2\|  \frac{x}{p} \|^2 ),
\end{equation}
where, in the last inequality, we used that fact that $1-y \le \exp(-y)$ for any $0 \le y \le 1$.

Consequently, we obtain the key inequality

\begin{equation} \label{eqn:fourier4}
\rho \le \frac{1}{p} \sum_{\xi \in \F_p} \prod_{i}|\cos \frac{ \pi v_i \xi}{p}  | \le  \frac{1}{p} \sum_{\xi \in F_p} \exp( - 2\sum_{i=1}^n  \|\frac{v_i \xi}{p}\|^2).
\end{equation}

{\it Large level sets.}  Now, we consider the level sets $S_m:=\{\xi| \sum_{i=1}^n  \| v_i \xi/p \| ^2 \le m  \} $.  We have

$$n^{-C} \le \rho  \le  \frac{1}{p} \sum_{\xi \in \F_p} \exp( -2 \sum_{i=1}^n  \| \frac{v_i \xi }{p} \| ^2) \le \frac{1}{p} + \frac{1}{p} \sum_{m \ge 1} \exp(-2(m-1)) |S_m| .$$

Because $\sum_{m\ge 1} \exp(-m) < 1$, there must be a large level set $S_m$ such that

\begin{equation}\label{eqn:level1} 
|S_m| \exp(-m+2) \ge  \rho  p. 
\end{equation}

In fact, because $\rho \ge n^{-C}$, we can assume that $m=O(\log n)$.

{\it Double counting and the triangle inequality.} By double -counting, we have

$$ \sum_{i=1}^n \sum_{\xi \in S_m} \|\frac{v_i \xi}{p} \| ^2 =   \sum_{\xi \in S_m} \sum_{i=1}^n  \|\frac{v_i \xi}{p} \| ^2 \le m |S_m |.$$

So, for most $v_i$

\begin{equation}\label{eqn:double1} 
\sum_{\xi \in S_m} \|\frac{v_i \xi}{p}\|^2  \le \frac{C_0 m }{n} |S_m| 
\end{equation}  

for some large constant $C_0$.

Set $C_0 = \eps^{-1}$. By averaging, the set of $v_i$ satisfying \eqref{eqn:double1} has a size of at least $(1-\eps)n$.  We call this set $V'$. The set $V\backslash V'$ has a size of at most $\eps n$, and this is the 
exceptional set that appears in Theorem \ref{theorem:optimal}. In the rest of the proof, we are going to show that $V'$ is a dense subset of a proper GAP.

Because $\|\cdot \|$ is a norm, by the triangle inequality, we have, for any $a  \in k V' $,

\begin{equation}\label{eqn:double2} 
\sum_{\xi \in S_m} \|\frac{a \xi}{p}\|^2  \le  k^2 \frac{C_0 m}{n} |S_m|. 
\end{equation}

More generally, for any $l  \le k $ and $a \in lV'$,

\begin{equation}\label{eqn:double3} 
\sum_{\xi \in S_m} \|\frac{a \xi}{p}\|^2  \le  k^2 \frac{C_0 m}{n} |S_m|. 
\end{equation}

{\it Dual sets.}  Define  $S_m^{\ast} :=\{ a | \sum_{\xi \in S_m}  \|\frac{a \xi}{p}\|^2 \le \frac{1}{200} |S_m |\}$ (the constant $200$ is ad hoc, and any sufficiently large constant would be sufficient).
$S_m^{\ast} $ can be viewed as some sort of a {\it dual} set of $S_m$. In fact, one can show, as far as cardinality is concerned, it does behave like a dual

\begin{equation}\label{eqn:dual1} 
|S_m^{\ast} |  \le \frac{8p}{|S_m |}. 
\end{equation}

To see this, define $T_a :=\sum_{\xi \in S_m} \cos \frac{2\pi a \xi}{p}$. Using the fact that $\cos 2\pi z \ge 1 -100 \|z\|^2 $ for any $z \in \R$, we have, for any $a \in S_m^{\ast}$

$$T_a \ge  \sum_{\xi \in S_m} (1- 100 \| \frac{a\xi}{p} \|^2 ) \ge \frac{1}{2} |S_m |. $$

However, using the basic identity $\sum_{a \in \F_p} \cos \frac{2\pi  ax}{p} = p\I_{x=0} $, we have

$$\sum_{a \in \F_p} T_a^2 \le 2p |S_m| . $$

\eqref{eqn:dual1} follows from the last two estimates and averaging.

Set $k := c_1 \sqrt{\frac{n}{m}}$, for a properly chosen constant $c_1= c_1(C_0) $. By \eqref{eqn:double3}, 
we have $\cup_{l=1}^k  l V'  \subset  S_m^{\ast} $. Set $V^{''} = V' \cup \{0\}$; we have 
$k V^{''} \subset S_m^{\ast} \cup \{0\} $. This results in the critical bound

\begin{equation}\label{eqn:dual2} 
|k V^{''} |  = O( \frac{p}{|S_m|}) = O(\rho^{-1} \exp(-m+2)).  
\end{equation}

{\it The long range inverse theorem.} The role of $\F_p$ is no longer important, so we can view the $v_i$ as integers.
The inequality \eqref{eqn:dual2} is exactly the assumption of the long range inverse theorem.

With this theorem in hand, we are ready to conclude the proof.
A slight technical problem is that $V^{''}$ is not a set but a multiset. Thus, we apply Theorem \ref{theorem:longrange} with 
$X$ as the set of distinct elements of $V^{''}$ (note that $kX=kV''$ if $k\ge 2$). Furthermore, $k = \Omega (\sqrt {\frac{n}{m}}) =\Omega(\sqrt{\frac{n}{\log n}})$, $\rho^{-1} \le n^{C}$ is bounded from above by $k^{2C+1}$.

It follows from Theorem \ref{theorem:longrange}  that $X$ is a subset of
a proper symmetric GAP $Q$ of rank $r=O_{C, \epsilon} (1)$ and cardinality

\begin{align*}
O_{C, \epsilon}(k^{-r}|kX|)= O_{C, \epsilon} (k^{-r} |kV^{''}|)&= O_{C, \epsilon} \left( \rho^{-1} \exp(-m) (\sqrt {\frac{n}{m}})^{-r}\right)\\
&= O_{C, \epsilon} (  \rho^{-1} n^{-r} ),
\end{align*}
concluding the proof.

\begin{remark} To prove Theorem \ref{theorem:optimal:general}, in the section describing {\it double counting and the triangle inequality}, we define $V'$ to be the collection of all $v_i\in V$ satisfying 

$$\sum_{\xi\in S_m} \|\frac{v_i\xi}{p}\|^2 \le \frac{m}{n'}|S_m|.$$ 

Next, with $k=c_1 \sqrt{\frac{n'}{m}}$ for some sufficiently small $c_1$, we obtain a bound similar to \eqref{eqn:dual2}, where $|kV''| = O(\rho^{-1} \exp(-m+2))$. We then conclude Theorem \ref{theorem:optimal:general} by applying the long range inverse theorem.  
\end{remark}

\section{Proof of Theorem \ref{theorem:continuous:NgV}}\label{section:continuous}

This proof will essentially follow the same steps as in the discrete case, with some additional simple arguments.

\noindent Given a real number $w$ and a variable $z$, we define the $z$-norm of $w$ by

$$\|w\|_z := (\E\|w(z_1-z_2)\|^2)^{1/2},$$

\noindent where $z_1,z_2$ are two iid copies of $z$. 

{\it Fourier analysis.}
Our first step is to obtain the following analogue of
 \eqref{eqn:fourier4}, using the Fourier transform.

\begin{lemma}[bounds for small ball probability]\label{lemma:upperboundforsmallball}

$$\rho_{r,z}(V)\le \exp(\pi r^2)\int_{\R^d}\exp(-\sum_{i=1}^n \|\langle v_i,\xi \rangle\|_z^2/2 - \pi \|\xi\|_2^2) \\d\xi .$$
\end{lemma}

This lemma is basically from \cite{TVcir}; the proof is presented in Appendix \ref{section:lemmasmallball}, for the reader's convenience.

Next, consider the multiset $V_{\beta}:=\beta^{-1}\cdot V =\{\beta^{-1}v_1,\dots, \beta^{-1}v_n\}.$ It is clear that

$$\rho_{\beta,z}(V) = \rho_{1,z}(V_\beta).$$

We now work with $V_\beta$. Thus $\rho_{1,z}(V_\beta)\ge n^{-O(1)}$ and $\sum_{v\in V_\beta}\|v\|^2 = \beta^{-2}$.

For concision, we write $\rho$ for $\rho_{1,z}(V_\beta)$. Set $M:= 2A \log n$, where $A$ is sufficiently large. From Lemma \ref{lemma:upperboundforsmallball} and the fact that $\rho\ge n^{-O(1)}$, we easily obtain

\begin{equation}\label{eqn:continuous:integral}
\int_{\|\xi\|_2\le M} \exp(-\frac{1}{2}\sum_{v\in V_\beta}\|\langle v,\xi \rangle \|_z^2-\pi \|\xi\|_2^2) \\d\xi\ge \frac{\rho}{2}.
\end{equation}

{\it Large level sets}. For each integer $0\le m \le M$, we define the level set

$$S_m:= \left \{\xi \in \R^d: \sum_{v\in V_\beta} \|\langle v,\xi \rangle \|_z^2 + \|\xi\|_2^2 \le m  \right \}.$$

Then, it follows from \eqref{eqn:continuous:integral} that $\sum_{m\le M} \mu(S_m) \exp(-\frac{m}{2}+1)\ge \rho$, where $\mu(.)$ denotes the Lebesgue measure of a measurable set. Hence, there exists $m\le M$ such that $\mu(S_m) \ge \rho\exp(\frac{m}{4}-2)$.

Next, because $S_m\subset B(0,\sqrt{m})$, by the pigeon-hole principle there exists a ball $B(x,\frac{1}{2})\subset B(0,\sqrt{m})$ such that

$$\mu (B(x,\frac{1}{2})\cap S_m ) \ge c_d\mu(S_m)m^{-d/2} \ge c_d\rho\exp(\frac{m}{4}-2)m^{-d/2}.$$

Consider $\xi_1,\xi_2\in B(x,1/2)\cap S_m$. By the Cauchy-Schwarz inequality (note that $\|.\|_z$ is a norm), we have

$$\sum_{v\in V_\beta} \|\langle v,(\xi_1-\xi_2) \rangle \|^2_z \le 4m.$$

Because $\xi_1-\xi_2 \in B(0,1)$ and $\mu(B(x,\frac{1}{2})\cap S_m - B(x,\frac{1}{2})\cap S_m) \ge \mu(B(x,\frac{1}{2})\cap S_m)$, if we put
$$T:=\{ \xi\in B(0,1), \sum_{i=1}^n \|\langle \xi,v_i\rangle \|^2_z \le 4m \},$$

\noindent then

$$\mu(T)\ge c_d\rho \exp(\frac{m}{4}-2)m^{-d/2}.$$

{\it Discretization}. Choose $N$ to be a sufficiently large prime (depending on the set $T$). Define the discrete box

$$B_0:=\left\{(k_1/N,\dots,k_d/N): k_i\in \Z, -N\le k_i \le N\right\}.$$

We consider all shifted boxes $x+B_0$, where $x\in [0,1/N]^d$. By the pigeon-hole principle, there exists $x_0$ such that the size of the discrete set $(x_0+B_0) \cap T$ is at least the expectation $|(x_0+B_0)\cap T| \ge N^d \mu(T)$ (to see this, we first consider the case when $T$ is a box).

Let us fix some $\xi_0\in (x_0+B_0)\cap T$. Then, for any $\xi \in (x_0+B_0)\cap T$, we have

$$\sum_{v\in V_\beta} \|\langle v,\xi_0-\xi\rangle \|^2_z \le 2 \left(\sum_{v\in V_\beta} \|\langle v,\xi\rangle \|^2_z + \sum_{v\in V_\beta} \|\langle v,\xi_0 \rangle \|^2_z\right) \le 16m.$$

Note that $\xi_0-\xi \in B_1:=B_0-B_0= \{(k_1/N,\dots,k_d/N): k_i\in \Z, -2N \le k_i \le 2N\}$. Thus, there exists a subset $S$ of size at least $c_dN^d\rho \exp(\frac{m}{4}-2)m^{-d/2}$ of $B_1$ such that the following holds for any $s\in S$:

$$\sum_{v\in V_\beta} \|\langle v,s\rangle \|^2_z \le 16m.$$

{\it Double counting}. We let $y=z_1-z_2$, where $z_1,z_2$ are iid copies of $z$. By the definition of $S$, we have

\begin{align*}
\sum_{s\in S} \sum_{v\in V_\beta} \|\langle v,s \rangle \|_z^2 &\le 16m |S|\\
\E_{y} \sum_{s\in S} \sum_{v\in V_\beta} \|y\langle v,s \rangle \|_{\R/\Z}^2 &\le 16m|S|.
\end{align*}

It is then implied that there exists $1\le |y_0|\le C_z$ such that

$$\sum_{s\in S} \sum_{v\in V_\beta} \|y_0\langle v,s \rangle \|_{\R/\Z}^2 \le 16m|S|\P(1\le |y|\le C_z)^{-1}.$$

However, by property \eqref{eqn:2ndmoment}, we have $\P(1\le |y|\le C_z)\ge 1/2$. Thus,

$$\sum_{s\in S} \sum_{v\in V_\beta} \|y_0\langle v,s \rangle \|_{\R/\Z}^2 \le 32m|S|.$$

 Let $n'$ be any number between $n^{\ep}$ and $n$. We say that $v\in V_\beta$ is {\it bad}  if

$$\sum_{s\in S} \|y_0 \langle v,s \rangle \|^2_{\R/\Z} \ge \frac{32m|S|}{n'}.$$

Then, the number of bad vectors is at most $n'$. Let $V_\beta'$ be the set of remaining vectors. Thus, $V_\beta'$ contains at least $n-n'$ elements.
In the remainder of the proof, we show that $V_\beta'$ is close to a GAP, as claimed in the theorem.

\vskip2mm

{\it Dual sets.} Consider an arbitrary $v\in V_\beta'$. We have $\sum_{s\in S} \|y_0 \langle s,v \rangle \|^2_{\R/\Z} \le 32m|S|/n'$.

Set $k:=\sqrt{\frac{n'}{64\pi^2 m}}$, and let $V_\beta'':=k(V_\beta'\cup \{0\})$. By the Cauchy-Schwarz inequality (see \eqref{eqn:double3}), for any $a\in V_\beta''$, we have

$$\sum_{s\in S} 2\pi^2 \|\langle s,y_0a \rangle \|^2_{\R/\Z} \le \frac{|S|}{2},$$

\noindent which implies

$$\sum_{s\in S}\cos(2\pi \langle s,y_0a\rangle) \ge \frac{|S|}{2}.$$

Observe that, for any $x\in C(0,\frac{1}{256d})$ (the ball of radius $1/256d$ in the $\|.\|_{\infty}$ norm) and any $s\in S\subset C(0,2)$, we always have $\cos(2\pi \langle s,x\rangle)\ge 1/2$ and $\sin(2\pi \langle s,x\rangle) \le 1/12$. Thus, for any $x\in C(0,\frac{1}{256d})$,

$$\sum_{s\in S}\cos\left(2\pi \langle s,(y_0a+x)\rangle\right) \ge \frac{|S|}{4}-\frac{|S|}{12} = \frac{|S|}{6}.$$

However,

\begin{align*}
\int_{x\in [0,N]^d} \left(\sum_{s\in S} \cos(2\pi \langle s, x\rangle )\right)^2 dx &\le \sum_{s_1,s_2\in S}\int_{x\in [0,N]^d} \exp \left( 2\pi \sqrt{-1}\langle s_1-s_2, x\rangle \right) dx\\
&\ll_d |S|N^d.
\end{align*}

Hence, we deduce the following:

$$\mu_{x\in [0,N]^d}\left ((\sum_{s\in S} \cos(2\pi \langle s, x\rangle ))^2 \ge (\frac{|S|}{6})^2 \right) \ll_d \frac{|S|N^d}{(|S|/6)^2} \ll_d \frac{N^d}{|S|}.$$

Now, using the facts that $S$ is large, $|S|\gg_dN^d \rho\exp(\frac{m}{4}-2)m^{-d/2}$ and $N$ was chosen to be large enough for $y_0V_\beta''+C(0,\frac{1}{256d}) \subset [0,N]^d$, we have

$$\mu ( y_0V_\beta''+C(0,\frac{1}{256d})) \ll_d \rho^{-1} \exp(-\frac{m}{4}+2)m^{d/2}.$$

Thus, we obtain the following analogue of (\ref{eqn:dual2}):

\begin{align}\label{eqn:sizeofkV'}
\mu \left( k(V_\beta'\cup \{0\})+C(0,\frac{1}{256dy_0}) \right) &\ll_d \rho^{-1} y_0^{-d}\exp(-\frac{m}{4}+2) m^{d/2}.
\end{align}

{\it The long range inverse theorem.}  Our analysis again relies on the long range inverse theorem. Let $D:=1024dy_0$. We approximate each vector $v'$ of $V_\beta'$ by its closest vector in $(\frac{\Z}{Dk})^d$,

$$\|v'-\frac{a}{Dk}\|_2 \le \frac{\sqrt{d}}{Dk}, \mbox{ with } a\in\Z^d.$$

Let $A_\beta$ be the collection of all such $a$. Because $\sum_{v'\in V_\beta'}\|v'\|^2_2 = O(\beta^{-2})$, we have

\begin{equation}\label{eqn:magnitudeV''}
\sum_{a\in A_\beta}\|a\|^2_2=O_{d,C_z}(k^2\beta^{-2}).
\end{equation}

It follows from \eqref{eqn:sizeofkV'} that

\begin{align*}
|k(A_\beta+C_0(0,1))|&=O_{d,C_z}\left(\rho^{-1}(Dk)^d y_0^{-d}\exp(-\frac{m}{4}+2)m^{d/2}\right)\\
&= O_{d,C_z}\left(\rho^{-1}k^d \exp(-\frac{m}{4}+2)m^{d/2}\right),
\end{align*}

\noindent where $C_0(0,r)$ is the discrete cube $\{(z_1,\dots,z_d)\in \Z^d: |z_i|\le r\}$.

Now, we apply Theorem \ref{theorem:longrange} to the set $A_\beta+C_0(0,1)$  (note that $0\in A_\beta$). That lemma implies there exists a proper GAP $P=\{\sum_{i=1}^r x_ig_i:|x_i|\le N_i\}\subset \Z^d$ containing $A_\beta+C_0(0,1)$ with a small rank $r=O(1)$ and small size

\begin{align*}
|P|&=O_{d,C_z}\left((\rho^{-1}k^d \exp(-\frac{m}{4}+2)m^{d/2} k^{-r}\right)\\
&=O_{d,C_z} (\rho^{-1}{n'}^{(-r+d)/2}).
\end{align*}

 Moreover, we learned from the proof of Theorem \ref{theorem:longrange} and Lemma \ref{lemma:Sarkozy} that $kQ$ can be contained in a set $ck(A_\beta+C_0(0,1))$ for some $c=O(1)$. Using \eqref{eqn:magnitudeV''}, we conclude that all generators $g_i$ of $Q$ are bounded,

$$\|g_i\|_2 =O_{d,C_z}(k \beta^{-1}).$$

Next, because $C_0(0,1)\subset Q$, the rank $r$ of $P$ is at least $d$. It is a routine calculation to see that $Q:=\frac{\beta}{Dk}\cdot P$ satisfies all of the 
required properties in Theorem \ref{theorem:continuous:NgV}.

\appendix

\section{Proof of the long range inverse theorem}\label{section:longrange}

The key lemma to prove our long range inverse theorem is an earlier result by Tao and the second author (\cite[Theorem 1.21]{TVjohn}), given below.

\begin{lemma}\label{lemma:4proper0}
Let $\ep>0,\gamma>0$ be constants. Assume that $X$ is a subset of integers such that $ |kX| \le k^\gamma|X|$ for some number
$k \ge 2$. Then, $kX$ is contained in a symmetric 2-proper GAP $Q$
with rank $r= O_{\gamma,\ep}(1)$ and cardinality $O_{\gamma,\ep}(|kX|)$.
\end{lemma}

Next, if $kX \subset kQ$, where $Q$ is a GAP, then it is natural to suspect that 
$X \subset Q$, but this is not always true. However, the conclusion holds if $kQ$ is 2-proper and $0 \in X$.

\begin{lemma}\label{lemma:4proper1} (Dividing sumsets relations)
Assume that $0 \in X$ and that $P=\{\sum_{i=1}^r
x_ia_i: |x_i|\le N_i\}$ is a symmetric 2-proper GAP that contains $kX$. Then $X\subset \{\sum_{i=1}^r x_ia_i: |x_i|\le
2N_i/k\}$.
\end{lemma}

A good way to keep this lemma in mind is the following. Consider the relation $X \subset P$. It is trivial that this 
relation can always be {\it multiplied}, namely, for all integers $k \ge 1$, $kX \subset kP$. The above lemma asserts that, under certain assumptions, the relation $kX \subset kP$ can be {\it divided}, giving $X \in P$.

\begin{proof}(of Lemma \ref{lemma:4proper1})
Without a loss of generality, we can assume that $k=2^l$. It is sufficient to 
show that $2^{l-1}X \subset \{\sum_{i=1}^r x_ia_i: |x_i|\le
N_i/2\}$. Because $0 \in X$, $2^{l-1}X \subset 2^lX\subset P$, any element
$x$ of $2^{l-1}X$ can be written as $x=\sum_{i=1}^r x_ia_i$, with
$|x_i|\le N_i$. Now, because $2x\in P \subset 2P$ and  $2P$ is proper (as $P$ is 2-proper), we must 
have $0\le |2x_i|\le N_i$.
\end{proof}

It is clear that Theorem \ref{theorem:longrange} follows from Lemma \ref{lemma:4proper0} and Lemma \ref{lemma:4proper1}.

\section{Remarks on  Theorem  \ref{theorem:continuous:NgV}} \label{section:corcont}

The purpose of this section is to give an example showing that the bound in Theorem \ref{theorem:continuous:NgV} cannot be improved and to provide a proof for Corollary \ref{corollary:continuous:NgV}.

First, consider the set $U := [-2n, -n] \cup [n, 2n]$. Sample $n$ points $v_1, \dots, v_n$ from $U$ independently with respect to the (continuous) uniform distribution, and let $A$ be the set of sampled points.
Let $\xi$ be the Gaussian random variable $N(0,1)$, and consider the sum

$$S:= v_1 \xi_1 + \dots + v_n \xi_n ,$$ where $\xi_i$ are iid copies of $\xi$.

$S$ has a Gaussian distribution with a mean $0$ and variance $\Theta (n^3)$, with a probability of one. Thus, for some interval $I$ of length $1$, 
$\P(S  \in I) \ge C n^{-3/2}$, for some constant $C$.

Set $n'= \delta n$, for some small positive constant $\delta$. Theorem \ref{theorem:continuous:NgV} states that (most of) $A$ is
$O(\frac{\log n}{\sqrt n})$-close to a GAP of rank $r$ and volume $O(n^{2 - \frac{r}{2}})$. We show that one cannot replace this bound by
$O(n^{2- \frac{r}{2}-\ep})$ for any $\ep$.  There are only three possible values for $r$: $r=1,2,3$. Our claim follows from the following simple lemma, whose proof remains as an exercise.

\begin{lemma} Let $C, \delta,\ep $ be positive constants and $n \rightarrow \infty$.
The following hold with a probability of $1-o(1)$ (with respect to the random choice of $A$).
\begin {itemize}

\item $A$ does not contain any subset of cardinality $(1-\delta)n$ that is $\frac{C\log n}{\sqrt{n}}$-close to a GAP of rank 1 and volume of at most $Cn^{3/2-\ep}$.

\vskip .1in

\item $A$ does not contain any subset of cardinality $(1-\delta)n$ that is  $\frac{C\log n}{\sqrt{n}}$-close to a GAP of rank  2 and volume of at most $Cn^{1-\ep}$.

\vskip .1in

\item  $A$ does not contain any subset of cardinality $(1-\delta)n$ that is  $\frac{C\log n}{\sqrt{n}}$-close to a GAP of rank 3 and volume of 
at most $Cn^{1/2-\ep}$.

\end{itemize}
\end{lemma}

\vskip2mm

The construction above can also be generalized to higher dimensions, but we do not attempt to do so here.

For the remainder of this section, we prove Corollary  \ref{corollary:continuous:NgV}. 

We consider the following two cases.

{\bf Case 1}: $r\ge d+1$. Consider the GAP $P$ at the end of the proof of Theorem
\ref{theorem:continuous:NgV}. Recall that $|P|=O_{d,C_z}(\rho^{-1}{n'}^{(d-r)/2})= O_{d,C_z}(\rho^{-1}/\sqrt{n'})$. Let

$$Q:=\frac{\beta}{Dk}\cdot P.$$

It is clear that $Q$ satisfies all of the conditions of Corollary \ref{corollary:continuous:NgV}. (Note that, in this case, we obtain a stronger approximation; almost all elements of $V$ are $O(\frac{\beta \log n'}{\sqrt{n'}})$-close to $Q$.)

{\bf Case 2}: $r=d$. Because the unit vectors $e_j=(0,\dots,1,\dots,0)$ belong to $P=\{\sum_{i=1}^d x_ig_i:|x_i|\le N_i\}\subset \Z^d$, the set of generators $g_i, i=1, \dots, d$ forms a base with the unit determinant of $\R^d$. In $P$, consider the set of lattice points with all coordinates divisible by $k$. We observe that (for instance, by \cite[Theorem 3.36]{TVbook}) this set can be contained in a GAP $P'$ of rank $d$ and cardinality $\max\left( O(\frac{1}{k^r}|P|,1\right)=\max \left(O(\rho^{-1}/{n'}^{r/2}),1 \right)$. (Here, we use the bound $|P|=O(\rho^{-1}\exp(-\frac{m}{4}) m^{d/2})$.) Next, define

$$Q:=\frac{\beta}{Dk}\cdot P'.$$

It is easy to verify that $Q$ satisfies all of the conditions of Corollary  \ref{corollary:continuous:NgV}. (Note that, in this case, we obtain a stronger bound on the size of $Q$.)

\section{Proof of  Lemma \ref{lemma:upperboundforsmallball} } \label{section:lemmasmallball}
We have

\begin{align*}
\P(\sum_{i=1}^n z_iv_i \in B(x,r)) &= \P(\|\sum_{i=1}^n z_iv_i -x\|_2^2 \le r^2)\\
&= \P\left(\exp(-\pi\|\sum_{i=1}^n z_iv_i-x\|_2^2) \ge \exp(-\pi r^2)\right)\\
&\le \exp(\pi r^2)\E \exp(-\pi\|\sum_{i=1}^n z_iv_i-x\|_2^2).
\end{align*}

Note that

$$\exp(-\pi\|x\|_2^2) = \int_{\R^d}e(\langle x,\xi \rangle) \exp(-\pi \|\xi\|_2^2) \\ d\xi.$$

We thus have

$$\P(\sum_{i=1}^n z_iv_i\in B(x,r))\le \exp(\pi r^2) \int_{\R^d}\E e(\langle \sum_{i=1}^n z_iv_i,\xi \rangle)e(-\langle x,\xi \rangle)\exp(-\pi \|\xi\|_2^2)\\ d\xi.$$

Using

$$|\E e(\langle \sum_{i=1}^n z_iv_i,\xi \rangle)| = \prod_{i=1}^n |\E e( z_i \langle v_i,\xi\rangle)|,$$

and

$$|\E e( z_i \langle v_i,\xi\rangle)| \le |\E e( z_i \langle v_i,\xi\rangle)|^2/2+1/2 \le \exp(-\|\langle v_i,\xi \rangle\|_z^2/2),$$

we obtain

$$\rho_{r,z}(V)\le \exp(\pi r^2) \int_{\R^d}\exp(-\sum_{i=1}^n\|\langle v_i,\xi \rangle \|_z^2/2 -\pi \|\xi\|_2^2) \\d\xi. $$

{\bf Acknowledgements.} The authors would like to thank K. Costello and the referees for carefully reading this manuscript and providing very helpful remarks.


\end{document}